# A Parametric Level Set Method for Topology Optimization based on Deep Neural Network (DNN)

Hao Deng, and Albert C. To[*]

Department of Mechanical Engineering and Materials Science, University of Pittsburgh, Pittsburgh, PA 15261

*Corresponding author. Email: albertto@pitt.edu

**Abstract**

This paper proposes a new parametric level set method for topology optimization based on Deep Neural Network (DNN). In this method, the fully connected deep neural network is incorporated into the conventional level set methods to construct an effective approach for structural topology optimization. The implicit function of level set is described by fully connected deep neural networks. A DNN-based level set optimization method is proposed, where the Hamilton-Jacobi partial differential equations (PDEs) are transformed into parametrized ordinary differential equations (ODEs). The zero-level set of implicit function is updated through updating the weights and biases of networks. The parametrized reinitialization is applied periodically to prevent the implicit function from being too steep or too flat in the vicinity of its zero-level set. The proposed method is implemented in the framework of minimum compliance, which is a well-known benchmark for topology optimization. In practice, designers desire to have multiple design options, where they can choose a better conceptual design base on their design experience. One of the major advantages of DNN-based level set method is its ability to generate diverse and competitive designs with different network architectures. Several numerical examples are presented to verify the effectiveness of the proposed DNN-based level set method.

## 1. Introduction

Topology optimization has experienced great attention and development in recent years in that it searches for the optimal material layout in a design domain with gradient-based algorithm. The first homogenization-based topology optimization method proposed by Bendsoe and Kikuchi is published in 1988 [1]. The research on topology optimization experienced a boost in the past two decades [2]. Solid Isotropic Material with Penalization (SIMP) method [3] is widely used in engineering because of effectiveness and simplicity. For SIMP method, the artificial density is applied to describe the material layout distribution, and optimal design is achieved through gradient-based optimization algorithm. However, the intermediate density may exist in optimal design, which blurs the boundaries and post-processing techniques are needed to remove gray area. In fact, it is hard for standard density-based method to eliminate intermediate density during optimization [4]. Based on standard density-based optimization framework, several advanced schemes are proposed in recent years to achieve feature control, robust design and length scale control [5-15], which is able to alleviate the aforementioned issue using various projection methods. Some other robust formulations [16-19] are also proposed in recent years to ensure manufacturability. In general, the post-processing techniques are needed to obtain a clear-boundary optimal design, while the performance of original optimal design may degenerate compared with the post-processed design.

Compared with SIMP method, the level set method is a "moving boundary" approach which evolves the boundaries of design in optimization with distinct boundaries [20]. In the beginning, Osher and Sethian [21] proposed the level set method to tackle the fronts of moving fluid. Level set methods represent the design boundaries using the zero-level set of an implicit function, and shape sensitivity analysis is implemented to compute the velocity field, which is incorporated into the Hamilton-Jacobi partial differential equation (PDE) to evolve the level set function. Osher and Sethian [22] proposed a level set method for design optimization problems based on the projected gradient approach. This work is further developed by Allaire [23] and Wang [24]. Allaire [23] proposed a numerical method based on a combination of shape derivative and level set method for front propagation, where the weight and perimeter constrains are considered as objective functions. Wang [24] used a scalar function of a higher dimension to represent the structural boundary with level set model, where a link between the velocity field and structural sensitivity analysis is identified. Compared with conventional level set method, several parametrized level set methods are proposed in recent years. Wang et al [25] incorporated radial basis functions (RBFs) into conventional level set method to construct a more efficient approach for topology optimization. An RBF-level set optimization method is proposed to transform the Hamilton-Jacobi PDE to parameterized ODEs, and level set function is updated through updating the expansion coefficients. Peng and Wang [26] proposed a piecewise constant level set (PCLS) method to resolve the shape and topology optimization problem, where the boundary is described by the PCLS functions. Jiang et al [27] applied a cardinal basis function to parameterize the level set function with unity collocation method, where a distance regularization energy functional is introduced to maintain the desired signed distance property in optimization. Recently, Guo et al [28, 29] proposed a new computational framework named moving morphable component (MMC), which embedded the moving morphable components into level set scheme. This computational scheme incorporates geometry and mechanical information into topology optimization in an explicit way and the structural complexity can be easily controlled in an explicit way. Another method called the Moving Morphable Voids (MMVs) was also proposed by Guo et al [30, 31], which introduced a set of geometry parameters to describe the boundary of the structure in an explicit way. In recent years, several other advanced parametrized level set method are proposed as described in Ref. [32]. Luo et al [33] proposed an efficient non-gradient approach for topology optimization without any sensitivity information. In his method, the material-field series expansion (MFSE) is applied to parametrize the geometry information, which achieves a considerable reduction of design variables, and the kriging-based optimization algorithm is implemented to resolve optimization problem based on the surrogate model.

Machine learning [34] has experienced a huge increase in research interest in the past decade, since it is a powerful tool for constructing the relationship between the input and output sampling data. With the dramatic growth of available data and development of new methods, machine learning has revolutionized our understanding of physical world such as image recognition [35], drug discovery [36], etc. Several successfully applications to physical problems can be found in recent years. Raissi et al [37] proposed physics-informed neural networks (PINNs) for solving partial differential equations. In his method, a physics-informed neural networks are trained to solve supervised learning tasks, which are constrained by given laws of physics. Based on this concept, several advanced machine learning-based method for solving PDEs [38-41] are proposed recently for forward and inverse PDE-based problems. Recent years have also witnessed several studies applying machine learning methods to resolve topology optimization problems. Yu et al [42] proposed a novel deep learning-based method to predict optimal design with a given boundary conditions without any iterative scheme. Lei et al [43] proposed a method to achieve real-time structural topology optimization based on machine learning method. His method combined MMC-based explicit framework with supported vector regression (SVR) to establish a mapping between design parameters and optimal designs. Oh and Jung et al [44] proposed an framework integrating topology optimization and generative models (generative adversarial networks) in an iterative way to generate new designs. As described in Ref. [45], deep neural networks are applied to represented shape's boundary with zero-level-set of the learned function, which can represent an entire class of shapes and thus the model size can be reduced by an order of magnitude compared with existing works. In the computational framework of parametrized level set method, a deep learning-based parameterized level set method is proposed in this paper to achieve topology optimization. The core of the current work is to incorporate deep neural networks into present level-set based topology optimization method. The implicit function is described by deep feed-forward neural networks. Thus, the sufficient smoothness and continuity of implicit function can be guaranteed.

At present, most topology optimization algorithm aim at finding optimal material layout which can minimize or maximize the objective function. In practice, the designers desire to generate multiple solutions that are diverse and competitive enough so that they can make a choice based on their experience from the aesthetic perspective or other functional requirements. Wang et al [46] achieve diverse and competitive designs by incorporating graphic diversity constraints, which is in the framework of SIMP method. Based on different penalty method, Yang et al [47] presents five simple and effective strategies to achieve multiple solutions, where these strategies are demonstrated to provide the designer with structurally efficient and topologically different solutions. Recently, He and Xie et al [48] proposed three stochastic approaches to generate diverse and competitive designs in the framework of bi-directional evolutionary structural optimization (BESO), where a series of random designs are produced with distinctly different topologies. For level set method, rare literatures are found in this field. Here, we proposed a DNN-based level set method to effectively generate diverse and competitive designs with high structural performance.

The paper is organized as follows. In Section 2, the implicit modeling based on deep neural network is presented. Section 3 describes the DNN level set topology optimization formulation in details. In Section 4, numerical examples are shown to illustrate the effectiveness of proposed parameterized level set method, followed by conclusions in Section 5.

## 2. Deep Neural Networks implicit modeling

To reconstruct the design domain with a single continuous and differentiable function, an implicit modeling method based on deep neural networks is presented here. The feedforward networks [49], with one or more layers between the input and output layers, are mainly used for function approximation. The typical architectures of DNNs are illustrated in Fig. 1, which contains input, hidden layers, and output. The mathematical formulation of deep feedforward neural networks can be defined as,

$$\mathbb{N}(x, y, \boldsymbol{\theta}) = \mathbb{N}\left(\boldsymbol{a}^{(L+1)}\left(\boldsymbol{h}^{(L)}\left(\boldsymbol{a}^{(L)}(\ldots \boldsymbol{h}^{(1)}(\boldsymbol{a}^{(1)}(x, y))\right)\right)\right) \quad (1)$$

where $\mathbb{N}$ denotes feedforward networks, and the $\boldsymbol{\theta}$ is parameter of network. The hidden layer is defined as $\boldsymbol{h}^{(l)}(\boldsymbol{x})$, a network with L hidden layers can be expressed as, where $\boldsymbol{a}^{(l)}(\boldsymbol{x})$ is a linear operation, expressed as,

$$\boldsymbol{a}^{(l)}(\boldsymbol{x}) = \boldsymbol{W}^{(l)}\boldsymbol{x} + \boldsymbol{b}^{(l)} \quad (2)$$

where $\boldsymbol{W}^{(l)}$ is weight matrix and $\boldsymbol{b}^{(l)}$ is bias vector for the $l\,th$ layer. The weight matrix $\boldsymbol{W}^{(l)}(l = 1, 2, \cdots L)$ and bias $\boldsymbol{b}^{(l)}(l = 1, 2, \cdots L)$ can be combined into a single parameter $\boldsymbol{\theta}$. $\boldsymbol{h}^{(l)}(l = 1, 2, \cdots L)$ are hidden-layer activation functions (kernel functions).

(a) One hidden layer          (b) Two hidden layer          (c) Three hidden layer

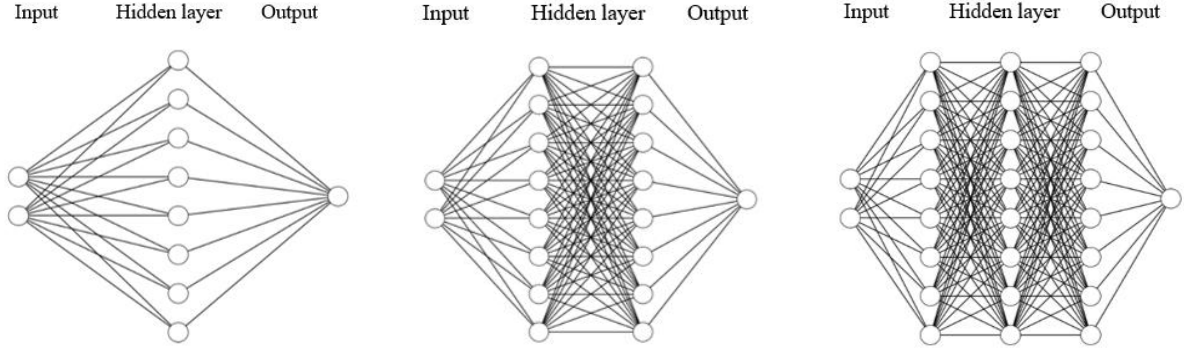

Figure 1. Architecture of DNNs

In fact, DNN is a universal approximator for nonlinear functions. It has been proven that a three-layer feedforward neural networks can approximate any continuous multivariate function to any accuracy [50]. DNN is a very effective tool for function approximation in high-dimensional spaces. Besides that, DNN models are analytically differentiable, and Graph-based automatic differentiation [51] techniques can be applied to easily gradient information. Compared with conventional discrete level set method, another extraordinary merit of DNNs is model reduction. As described by literature, DeepSDF [45] using deep neural networks to learn signed distance function to represent the zero-level set of complex geometry shape demonstrates extraordinary model reduction ability, where the high fidelity of shape using DeepSDF is validated. The initial weights of a networks are critical for convergence of a training. In practice, it is common to initialize all the weights and biases with random zero-mean values [52]. The initialization of an implicit level set function represented by DNNs can be formulated as follows,

$$\begin{cases} Find: \boldsymbol{\theta} \\ Min: \sum_{i=1}^{N} \|\mathbb{N}(x_i, y_i, \boldsymbol{\theta}) - \Phi(x_i, y_i)\|_2 \end{cases} \quad (3)$$

where $\mathbb{N}$ is the feedforward neural network, and $\Phi$ is the target implicit function. Operator $\|\cdot\|_2$ denotes 2-norm. $(x, y)$ denotes the coordinate of point. The backpropagation learning algorithm [52] is applied here to train the neural networks. The activation function is chosen as hyperbolic tangent function. Each layer of network contains 8 neurons. The target implicit function is plotted in Fig. 2. The training results of the plate with 5 circle holes inside using three different architectures are presented in Fig. 3. Note that for network with one hidden layer, the trained shape cannot achieve high fidelity as shown in Fig. 3(a), while the better shape training result can be found in Fig. 3(c) for network with three hidden layers.

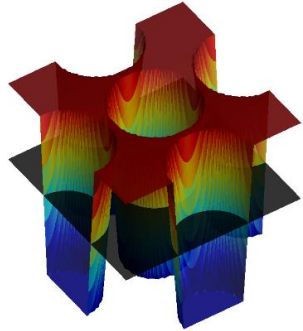

Figure 2. Target implicit function

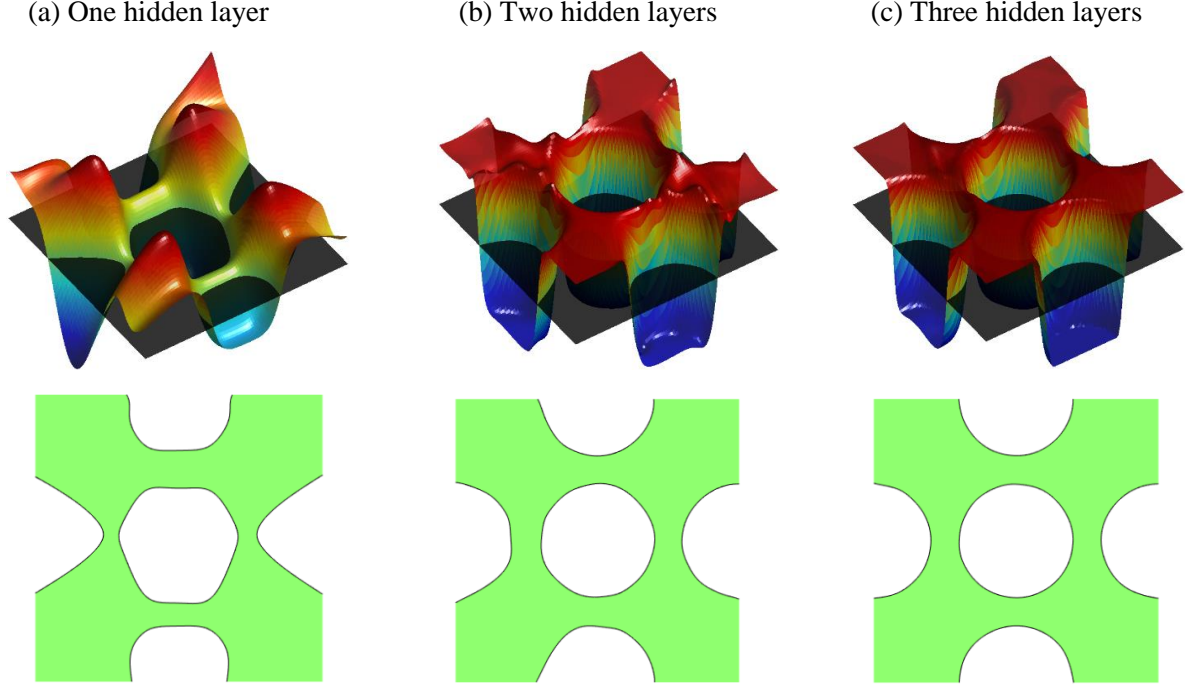

Figure 3. Training results (a) one hidden layer (b) two hidden layers (c) three hidden layers

## 3. DNN level set method for structural topology optimization

### 3.1 Conventional level set-based topology optimization

Conventional level set method uses a zero contour (2D) or isosurface (3D) to represent the boundaries of geometry, which is introduced by Osher and Sethian [21] to simulate the motion of dynamic interfaces. The interface is described by the zero-level set of implicit function $\Phi(x)$, which is Lipschitz-continuous in design domain. In this paper, the level set function $\Phi(x)$ is defined as:

$$\begin{cases} \Phi(\pmb{x},t) > 0, & (\pmb{x} \in \Omega) \\ \Phi(\pmb{x},t) = 0, & (\pmb{x} \in \partial\Omega) \\ \Phi(\pmb{x},t) < 0, & (\pmb{x} \in D\backslash\Omega) \end{cases} \quad (4)$$

where $D$ is design domain, $\Omega$ represents all admissible shapes, $\partial\Omega$ denotes the boundary of shape, and $t$ is the pseudo time [23] of shape dynamic evolution. The Hamilton-Jacobi PDE can be obtained through differentiating the zero-level set with respect to pseudo time $t$ as follows:

$$\frac{\partial \Phi}{\partial t} - V_n |\nabla \Phi| = 0 \quad (5)$$

where $V_n$ is normal velocity computed through sensitivity analysis. The shape of zero-level set evolves along the gradient direction through solving the above Hamilton-Jacobi equation. This equation can be resolved using upwind schemes, where a reinitialization procedure is needed as an auxiliary step to avoid the implicit function becomes too flat or steep in the vicinity of its zero-level set. In this paper, the objective is chosen as minimizing structural compliance $J(\Phi)$, which can be formulated as follows:

$$min: J(\Phi) = \int_D \big(\pmb{\varepsilon}(\pmb{u})\colon \pmb{C}\colon \pmb{\varepsilon}(\pmb{u})\big) H(\Phi) d\Omega \quad (6)$$

$$\text{Subject to} \begin{cases} Volume = \int_D H(\Phi)d\Omega = V_0 \\ a(\boldsymbol{u},\boldsymbol{v},\Phi) = l(\boldsymbol{v},\Phi) \\ \boldsymbol{u} = \boldsymbol{u}_0 \quad in\ \Gamma_u \\ \boldsymbol{C}:\boldsymbol{\varepsilon}(\boldsymbol{u})\cdot\boldsymbol{n} = \boldsymbol{\tau} \quad in\ \Gamma_\tau \end{cases} \tag{7}$$

where the notations in above equations can be written as,

$$a(\boldsymbol{u},\boldsymbol{v},\Phi) = \int_D \left(\boldsymbol{\varepsilon}(\boldsymbol{u}):\boldsymbol{C}:\boldsymbol{\varepsilon}(\boldsymbol{v})\right) H(\Phi)d\Omega \tag{8}$$

$$l(\boldsymbol{v},\Phi) = \int_{\Gamma_\tau} \boldsymbol{\tau}\cdot\boldsymbol{v}d\Gamma + \int_D \boldsymbol{b}\cdot\boldsymbol{v}H(\Phi)\,d\Omega \tag{9}$$

where $\boldsymbol{u}$ is the displacement and $\boldsymbol{\varepsilon}(\boldsymbol{u})$ is the strain. $H(\cdot)$ denotes the Heaviside step function. $\boldsymbol{C}$ is the elastic tensor. $\Gamma_u$ and $\Gamma_\tau$ denote the displacement and force boundary, respectively. $\boldsymbol{u}_0$ is prescribed displacement boundary conditions. $\boldsymbol{\tau}$ and $\boldsymbol{b}$ denote the traction at boundary and body force in the domain, respectively. $a(\boldsymbol{u},\boldsymbol{v},\Phi)$ is the energy bilinear form and $l(v,\Phi)$ is the load linear form. $\boldsymbol{v}$ is a virtual displacement field. The operator $(:)$ denotes tensor contraction. For Heaviside step function, which equals zero in the void area and one in the solid area. In practice, the step Heaviside function is approximated by a smooth function to ensure the differentiable in the transition area. The smoothed Heaviside function can be formulated as follows,

$$H_\delta(\Phi) = \begin{cases} \delta & \Phi < -\Delta \\ 0.75(1-\delta)\left(\frac{\Phi}{\Delta} - \frac{\Phi^3}{3\Delta^3}\right) + \frac{1+\delta}{2} & -\Delta \leq \Phi \leq \Delta \\ 1 & \Phi > \Delta \end{cases} \tag{10}$$

The $\delta$ is a small value, and $\Delta$ denotes the half of the transition width. Detailed mathematical properties of smoothed Heaviside function can be found in Ref. [24]

### 3.2 DNN level set optimization method

A DNN level set method is proposed to convert the Hamilton-Jacobi PDE into system of ODEs in the design domain for topology optimization. For conventional level-set method, the implicit function $\Phi(x)$ is updated through solving Hamilton-Jacobi equations to obtain the optimal topology. In this paper, DNN implicit modeling is applied to represent implicit function $\Phi(x)$, where the evolution of $\Phi(x)$ is equivalent to updating the parameters of networks. As mentioned in previous section, the implicit function represented by the time-dependent neural networks can be expressed as

$$\Phi(x,t) = \mathbb{N}(x,y,\boldsymbol{\theta}(t)) \tag{11}$$

Substituting Eq. (11) into the Hamilton-Jacobi equation, the parametrized ODE can be written as,

$$\frac{\partial \mathbb{N}(x,y,\boldsymbol{\theta}(t))}{\partial t} - V_n\left|\nabla\big(\mathbb{N}(x,y,\boldsymbol{\theta}(t))\big)\right| = 0 \tag{12}$$

$$\frac{\partial \mathbb{N}(x,y,\boldsymbol{\theta})}{\partial \boldsymbol{\theta}}\frac{\partial \boldsymbol{\theta}(t)}{\partial t} - V_n\left(\left(\frac{\partial\big(\mathbb{N}(x,y,\boldsymbol{\theta}(t))\big)}{\partial x}\right)^2 + \left(\frac{\partial\big(\mathbb{N}(x,y,\boldsymbol{\theta}(t))\big)}{\partial y}\right)^2\right)^{1/2} = 0 \tag{13}$$

the Moore–Penrose inverse $\mathcal{M}^+$ of Matrix $\frac{\partial \mathbb{N}(x,y,\boldsymbol{\theta})}{\partial \boldsymbol{\theta}}$ is applied to obtain the least squares solution of above system

$$\frac{\partial \boldsymbol{\theta}(t)}{\partial t} - \mathcal{M}^+ V_n \left( \left(\frac{\partial(\mathbb{N}(x,y,\boldsymbol{\theta}(t)))}{\partial x}\right)^2 + \left(\frac{\partial(\mathbb{N}(x,y,\boldsymbol{\theta}(t)))}{\partial y}\right)^2 \right)^{1/2} = 0 \quad (14)$$

where the $\mathcal{M}^+$ can be expressed as,

$$\mathcal{M}^+ = \left( \left(\frac{\partial \mathbb{N}(x,y,\boldsymbol{\theta})}{\partial \boldsymbol{\theta}}\right)^T \frac{\partial \mathbb{N}(x,y,\boldsymbol{\theta})}{\partial \boldsymbol{\theta}} \right)^{-1} \left(\frac{\partial \mathbb{N}(x,y,\boldsymbol{\theta})}{\partial \boldsymbol{\theta}}\right)^T \quad (15)$$

In Equation (14), the DNN coefficients are time dependent, where the initial value of DNN can be obtained through BP training. Therefore, a PDE problem is transformed into the ODE problem. It is worth to mention that the derivative information of DNN with respect to its parameter or input is ready to be obtained through graph-based automatic differentiation. To obtain high accuracy and stable solutions, Equation (14) is solved by a Runge-Kutta-Fehlberg (RKF45) method [53], which is recommended by Ref [54]. RKF45 method is able to determine if the proper step size is being used, where two different approximations for the solution are made and compared. If the two approximations do not agree to a specified accuracy, the step size is reduced. More details of this method can be found in Ref [53]. In general, the timestep size should be sufficiently small to achieve the numerical stability due to the Courant–Friedrichs–Lewy (CFL) condition for stability [55].

Based on shape derivative, the normal velocity $V_n$ along the free moving boundary form minimum compliance can be expressed as follows:

$$V_n = \boldsymbol{\varepsilon}(\boldsymbol{u}) : \boldsymbol{C} : \boldsymbol{\varepsilon}(\boldsymbol{u}) - \lambda \quad (16)$$

where $\lambda$ is the Lagrange multiplier to enforce the constraint of volume fraction. The augmented Lagrangian updating scheme [56] is applied to update $\lambda$ in this paper.

### 3.3 Parametrized reinitialization scheme

For standard level set method, irregularities during evolution may occur, which inevitable leads to instability of level set evolution. To overcome this difficulty, reinitialization scheme was introduced to regularize LSF to maintain the stability of boundary evolution. In general, reinitialization is implemented by periodically reshaping LSF as a signed distance function [57]. A standard method for initialization is obtaining the steady solution of the following evolution equation:

$$\frac{\partial \Phi}{\partial t} = sign(\Phi_{initial})(1 - |\nabla \Phi|) \quad (18)$$

$$\Phi_{initial} = \Phi(t = 0)$$

where $\Phi_{initial}$ is the level set function to be initialized, and $sign(\cdot)$ denotes the sign function. Substituting Eq. (11) into reinitialization in Eq. (18), the following equation can be obtained:

$$\frac{\partial \mathbb{N}(x,y,\boldsymbol{\theta})}{\partial \boldsymbol{\theta}} \frac{\partial \boldsymbol{\theta}}{\partial t} - sign(\Phi_{initial})(|\nabla \mathbb{N}(x,y,\boldsymbol{\theta})| - 1) = \mathbf{0} \quad (19)$$

The reinitialization procedure used in Eq. (19) usually slightly moves the zero-level set contour, which may cause inconsistencies during optimization process. In general, the reinitialization procedure needs to be implemented periodically to maintain the signed-distance function [13, 58]. Hartmann et al [59] proposed a constrained reinitialization scheme, where the least-squares solution is obtained to keep the location of zero-level contour. Long et al [32] proposed a new level set scheme using double-well potential function to achieve distance regularization inside the topology optimization loop. In this paper, the regular reinitialization scheme is applied periodically to avoid the implicit function becoming too flat or too steep. More details of implementation can be found in Ref [23]. Equation (19) can be solved using the Runge-Kutta-Fehlberg (RKF45) method [55]. The proposed DNN-based level set algorithm is described as follows,

1. Initialization of the parametrized level-set function $\Phi(\boldsymbol{\theta})$, corresponding to an initial guess $\Phi_0$. The initial weights and biases are determined through back-propagation algorithm (Eq. (3)).
2. Iteration until convergence, for $k \geq 0$
   (a) Based on implicit function $\Phi_k$, the material distribution in the design domain can be computed through Heaviside function Eq. (10). Solve the equilibrium equation (Eq. (7)) using FEM technique to obtain displacement field.
   (b) Compute objective $J(\Phi)$ and the normal velocity $V_n$. Updating the weights and biases by solving the parametrized Hamilton-Jacobi equation Eq. (13) using RBF45 scheme.
   (c) For stability reasons, the reinitialization of the level set function $\Phi$ by solving Eq. (19) in every iteration.

The flowchart of algorithm can be found in Fig. 4.

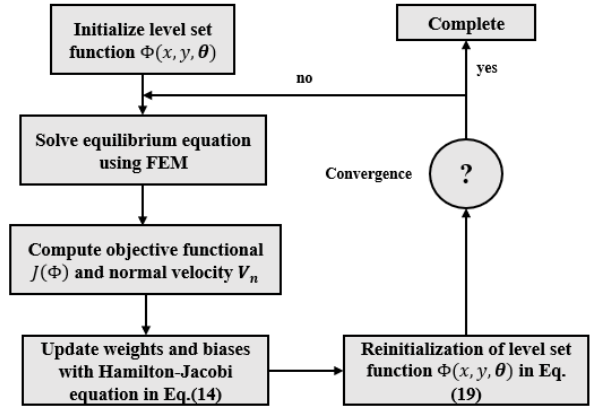

Figure 4. Flowchart of level set method based on DNN

## 4. Numerical examples

In this section, several two-dimensional numerical examples are demonstrated to verify the effectiveness of the proposed DNN level set method. Unless stated otherwise, the following parameters are chosen as: the elastic modulus $E = 1$ for solid materials, modulus $E = 1 \times 10^{-6}$ for void materials. The Poisson's ratio is chosen as $v = 0.3$. The volume fraction constraint is set to be 0.4. The deep neural networks are initialized using backpropagation training algorithm. The activation function is chosen as hyperbolic tangent function [60]. The detailed initialization description is presented in section 2. For all numerical examples, the design domain is discretized by rectilinear mesh with grid size equals 1. The finite element analysis is based on an 'ersatz material' approach, which is a well-known method for level set topology optimization [23].

### 4.1 MBB beam

The MBB beam is investigated in the present example for minimum compliance problem. The boundary condition is plotted in Fig. 5, where a concentrated force $P = 1$ is applied at the mid of the top edge. The architecture of network is presented in Fig. 6. Note that a fixed support is located at bottom-left corner, while a roller is at the bottom-right corner. The design domain is meshed with $200 \times 100$ elements with a grid size equals to 1. The fixed Lagrange multiplier is chosen as $l = 5$, and the time step is chosen as $\tau = 3 \times 10^{-3}$. For the first case, the implicit function is represented by the neural network with one hidden layer, where every hidden layer contains 8 neurons (Fig. 6). The total number of design variables is 33. To generate a symmetry design, only half of design domain is used for optimization, where symmetric boundary condition is implemented. The final optimal design is shown in Fig. 8(a). The initial design

(training result) is plotted in Fig. 7(a). Note that the neural network with one hidden layer is a shallow network, and the hole in training result is not a perfect circle due to limited fitting ability. A stable topology optimized solution can be achieved through solving the ODEs in Eq. (14) after 150 iterations (Fig.10(a)).

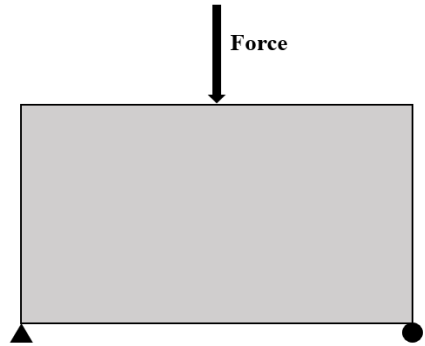

Figure 5. Compliance design of an MBB beam

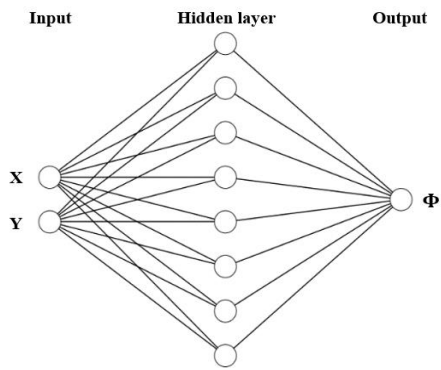

Figure 6. Architecture of networks

(a) hidden layers: 8  (b) hidden layers: $8 \times 8$  (c) hidden layers: $8 \times 8 \times 8$

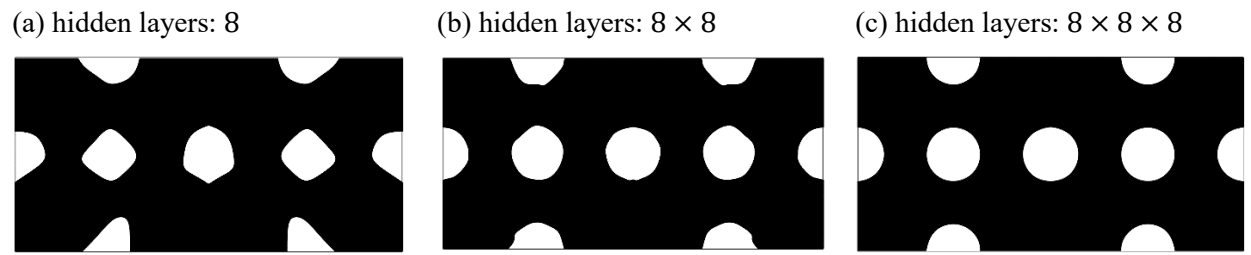

Figure 7. Initial design of an MBB beam design (DNN with one hidden layer)

(a) Comp:33.39 (8)  (b) Comp:32.57 ($8 \times 8$)  (c) Comp: 36.02 ($8 \times 8 \times 8$)

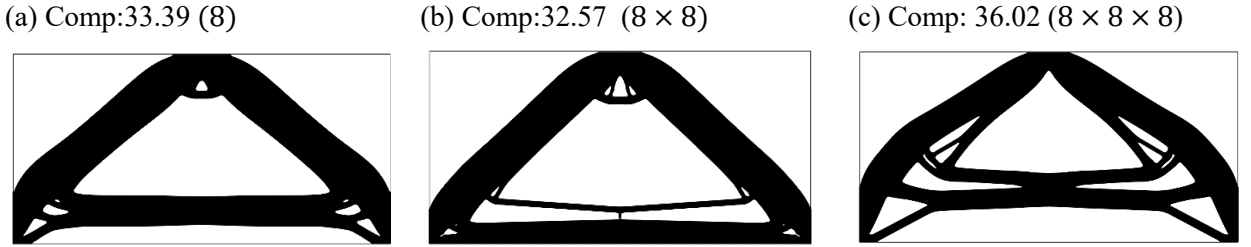

Figure 8. Optimized design of an MBB beam

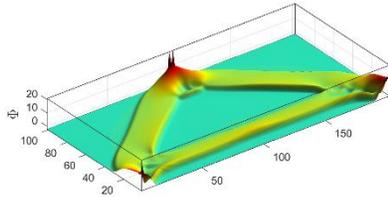 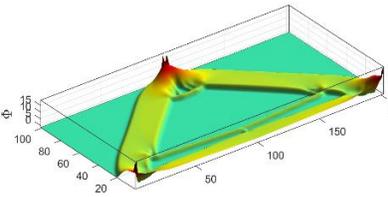 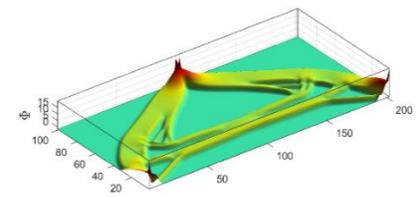

(a) hidden layers: 8      (b) hidden layers: 8 × 8      (c) hidden layers: 8 × 8 × 8

Figure 9. Implicit function of optimized design

To make a comparison, neural networks with two hidden layers and three hidden layers are chosen to represent the implicit function of level set. The optimization parameter setting is the same as the previous example. The architecture of network is shown in Fig. 1(b) and (c). The inputs are the point coordinates $(x, y)$ in the design domain, and the output is the value of the implicit function at the present point. The initial training results are presented in Fig. 7(b) and 7(c), and the optimal designs are displayed in Fig. 8(b) and 8(c). The total number of design variables are 105 (two hidden layers) and 177 (three hidden layers). The implicit function of optimal design is shown in Fig. 9. The compliance values of two optimal designs are 32.57 and 39.42, respectively. Convergence history is presented in Fig. 10. Since the level set function is updated via updating the parameters of neural networks, new holes can be generated freely from the mathematical point of view. This salient point is verified in Fig. 8, where new small holes are generated during the optimization process. A benchmark design generated by SIMP method using top88 code as provided in Ref. [61] is presented in Fig. 11. Note that the radius of filter for standard SIMP method is chosen as $r = 2$. The compliance of optimal design produced by SIMP method is 33.48. Compared with solutions from DNN-based level set method, the structural compliance values of optimized designs (networks with one, two layers) are quite close with the benchmark solution.

(a) hidden layers: 8

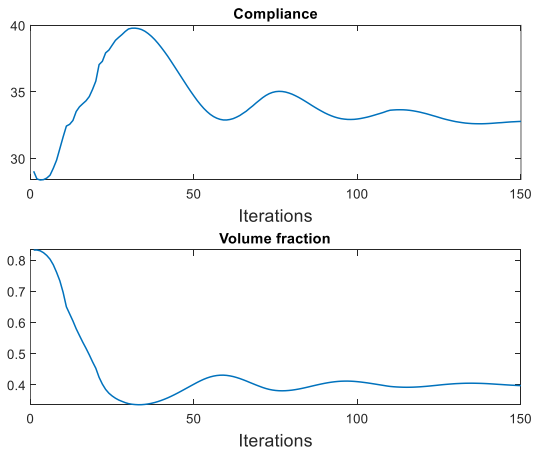

(b) hidden layers: $8 \times 8$

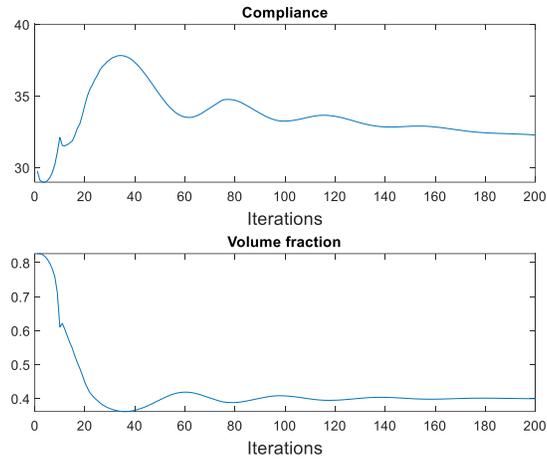

(c) hidden layers: $8 \times 8 \times 8$

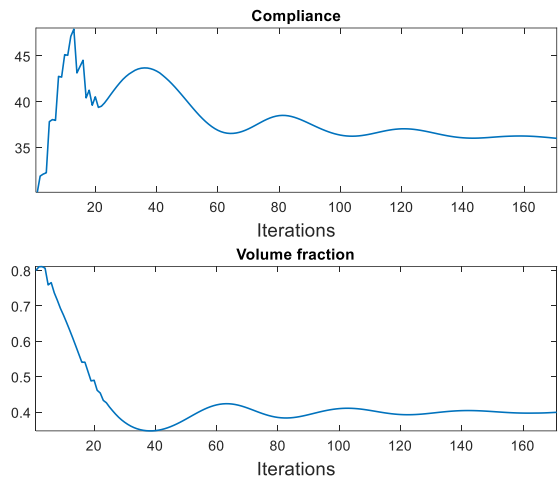

Figure 10. Convergence history (a) hidden layers: 8 (b) hidden layers: $8 \times 8$ (c) hidden layers: $8 \times 8 \times 8$

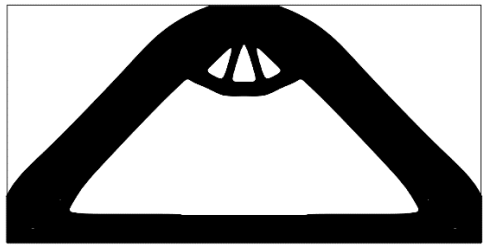

Figure 11. Benchmark design (Compliance: 33.48)

To further verify the effectiveness of proposed DNN-based level set method for diverse and competitive designs. Different network architectures are selected to produce diverse designs as shown in Fig. 12. It is worth mentioning that different network is capable of generating distinctly different solutions. For the single layer network with 5 neurons, the optimized result is simple and no intricate geometric details are found as shown in Fig. 5(a). For network with 2 layers (15 × 15), the optimized design seems to be more complex with several truss-like support components inside as shown in Fig. 12(f). For a given architecture of neural network (NN), the implicit function (segmentation) represented by NN is a subspace of all possible designs. Thus, NN with different architectures describes different subspace of solutions, which give the explanation of diverse solutions.

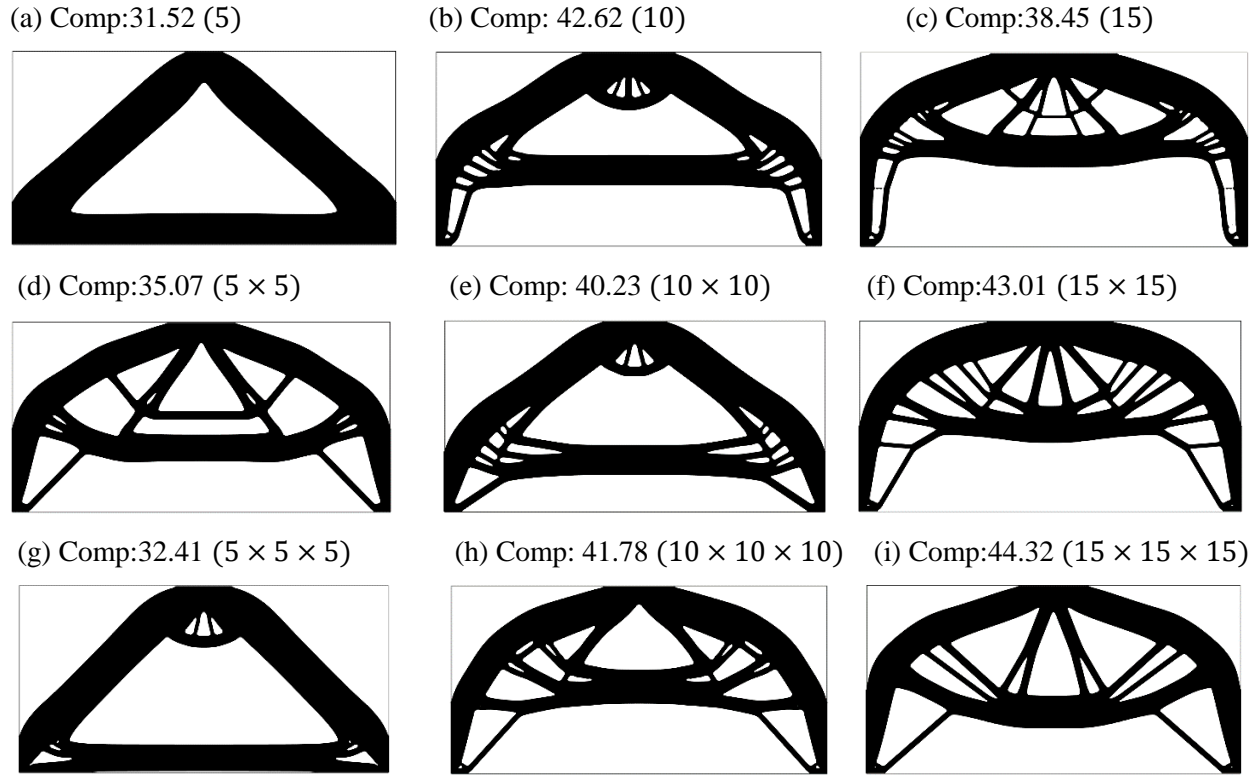

(a) Comp:31.52 (5)  (b) Comp: 42.62 (10)  (c) Comp:38.45 (15)
(d) Comp:35.07 (5 × 5)  (e) Comp: 40.23 (10 × 10)  (f) Comp:43.01 (15 × 15)
(g) Comp:32.41 (5 × 5 × 5)  (h) Comp: 41.78 (10 × 10 × 10)  (i) Comp:44.32 (15 × 15 × 15)

Figure 12. Diverse and competitive designs generated by DNN-based level set method

## 4.2 Short cantilever beam

The minimum compliance design of a short cantilever beam is presented in Fig. 13. The design domain is a square with fixed boundary condition on the left side, and a vertical concentration force $F = 1$ is applied at the midpoint of right side. The design domain is meshed with $100 \times 100$ elements with grid size equals 1. The optimization formulation is described in Eq. (6) and (7). A fixed Lagrange multiplier $l = 3$ is applied for volume constraint, and the time step is chosen as $\tau = 3 \times 10^{-3}$. For the first case, neural network with one hidden layer is chosen to represent the implicit function, where the hidden layer contains 8 neurons. Because of the limited fitting ability of shallow neural network, the initialized shape has some artifacts near the boundary, and the training result is presented in Fig. 14(a). The final optimal design is displayed in Fig. 15(a). The implicit function of optimal design is shown in Fig. 16(a), and the optimization progress converges after 120 iterations (Fig. 17(a)). To make a comparison with numerical results generated by shallow neural networks, the networks with two or three hidden layers are examined here. The architectures of networks are shown in Fig. 1(b) and Fig. 1(c), where each layer contains 8 neurons. The total number of design variables is 105 (two hidden layers) and 177 (three hidden layers). The other optimization settings are the same as before. The optimal layout and implicit functions are demonstrated in Fig. 15(b-c) and Fig. 16(b-c). Obviously, the optimal designs using deep neural networks (two or three hidden layers) have more intricate geometric features compared to the result obtained by shallow neural network. The convergence histories for the two designs are demonstrated in Fig. 17(b) and (c), where the stable topology optimized design is achieved after 140 iterations. A benchmark design obtained with the SIMP method is demonstrated in Fig. 18. It is worth to mention that designs produced by the NN with 2 or 3 layers is slightly lower than the value of benchmark design (around 3.5% difference).

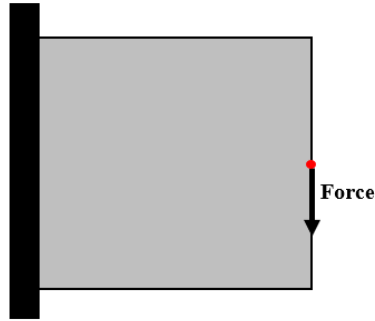

Figure 13. Compliance design of an short cantilever beam

(a) hidden layers: 8          (b) hidden layers: 8 × 8          (c) hidden layers: 8 × 8 × 8

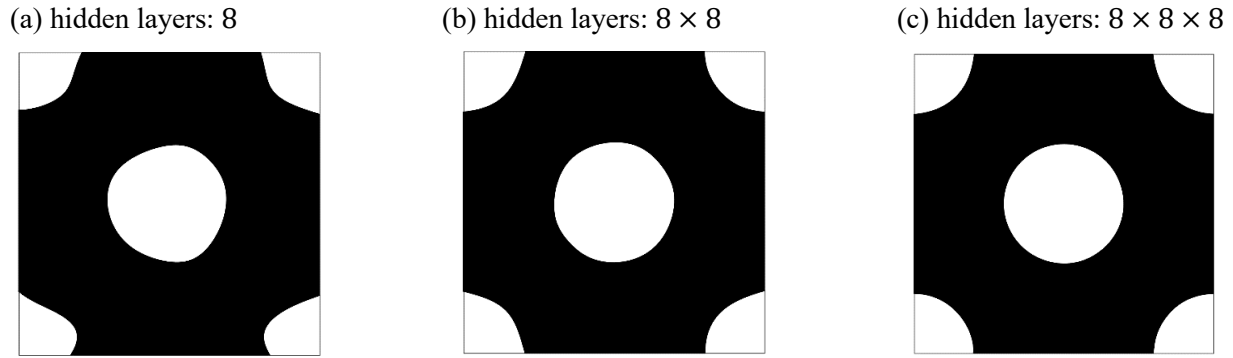

Figure 14. Initial design of an short cantilever beam

(a) hidden layers: 8 (181946)    (b) hidden layers: 8 × 8 (176893)  (c) hidden layers: 8 × 8 × 8 (177848)

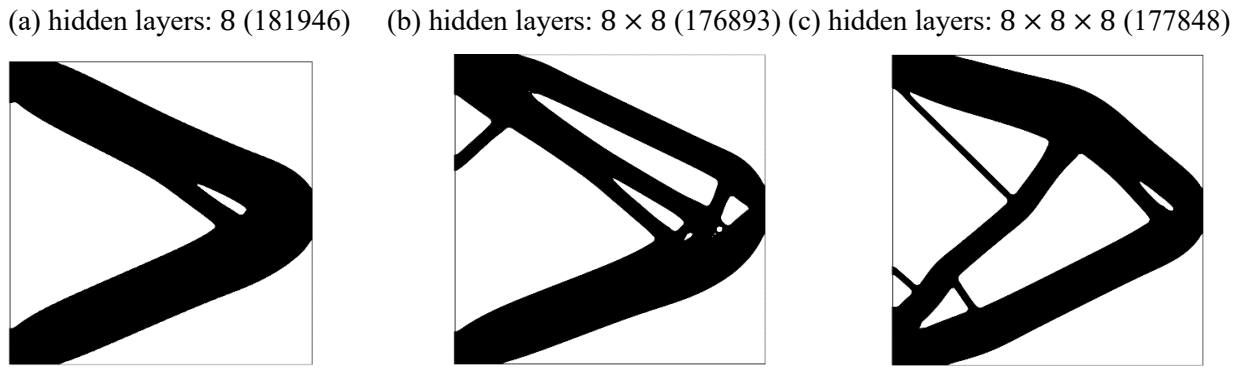

Figure 15. Optimized design of an short cantilever beam

| (a) hidden layers: 8 | (b) hidden layers: 8 × 8 | (c) hidden layers: 8 × 8 × 8 |
|---|---|---|
| 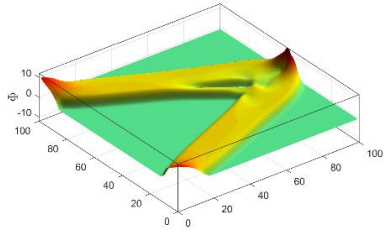 | 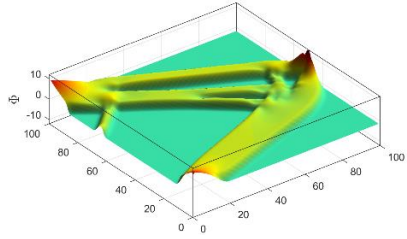 | 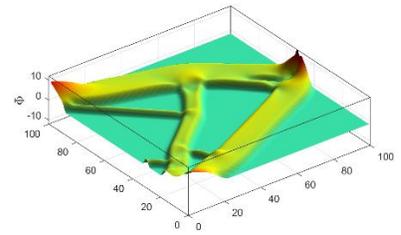 |

Figure 16. Implicit function of optimized design

(a) hidden layers: 8

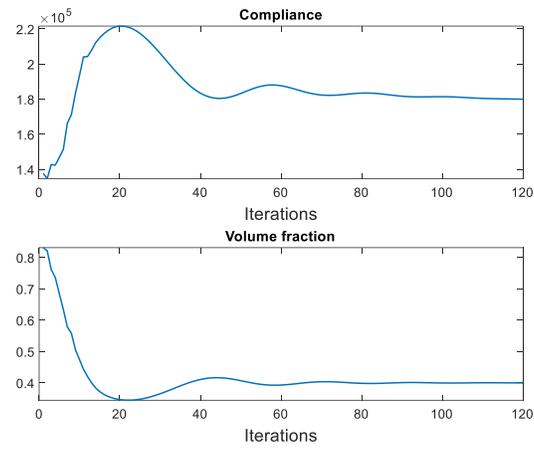

(b) hidden layers: 8 × 8

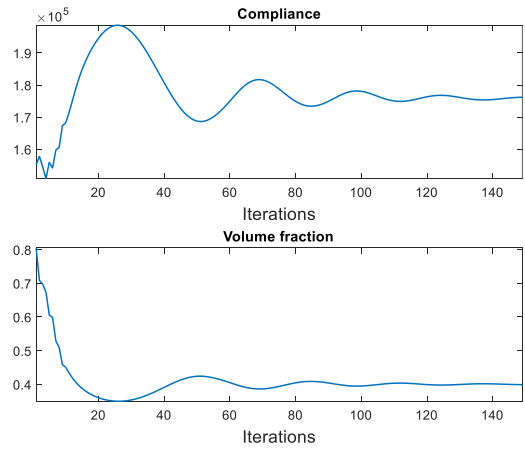

(c) hidden layers: 8 × 8 × 8

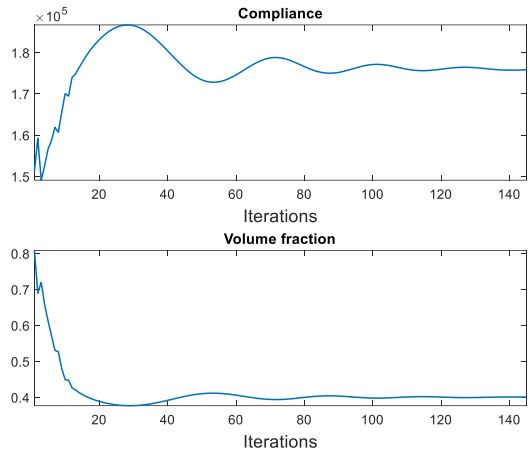

Figure 17. Convergence history (a) hidden layers: 8  (b) hidden layers: 8 × 8  (c) hidden layers: 8 × 8 × 8

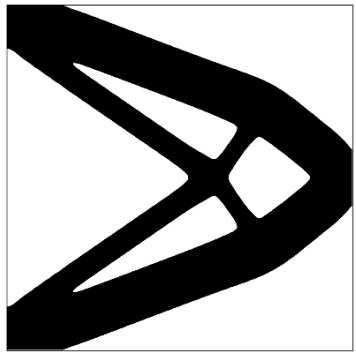

Figure 18. Benchmark design (Compliance: 183543)

To further generate multiple alternatives, 9 different NN architectures are examined to obtain the solutions as shown in Fig. 19. Obviously, although these designs have distinctly different topologies, the compliance value of optimized design are close with respect to benchmark (difference less than 5%). Unsymmetrical designs can be easily obtained as plotted in Fig.19 (b)-(i).

(a) Comp:174319 (5)          (b) Comp: 176569 (10)          (c) Comp:175819 (15)

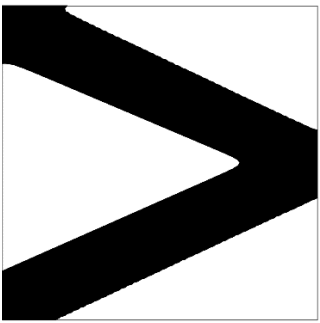 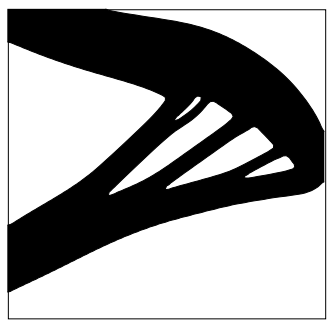 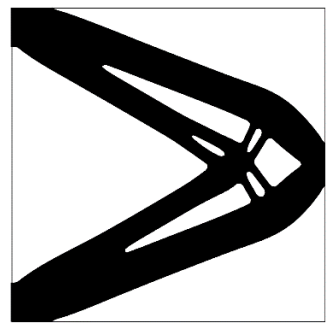

| (d) Comp:197768 (5 × 5) | (e) Comp: 174970 (10 × 10) | (f) Comp:181101 (15 × 15) |
|---|---|---|
| 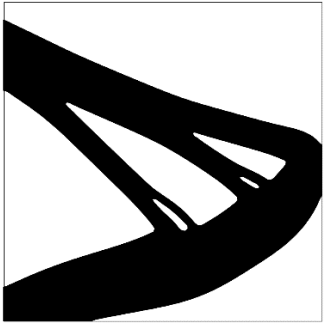 | 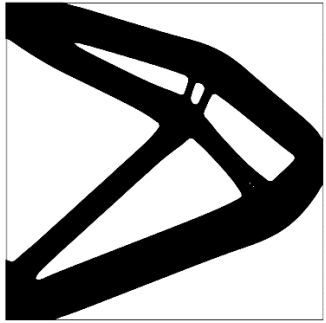 | 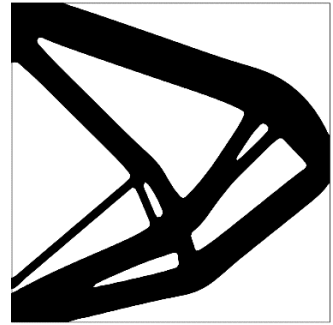 |
| (g) Comp:190953 (5 × 5 × 5) | (h) Comp: 176683 (10 × 10 × 10) | (i) Comp:178011 (15 × 15 × 15) |
| 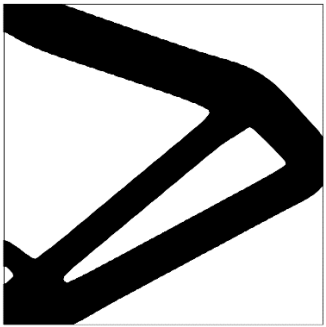 | 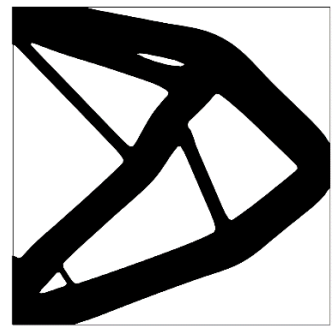 | 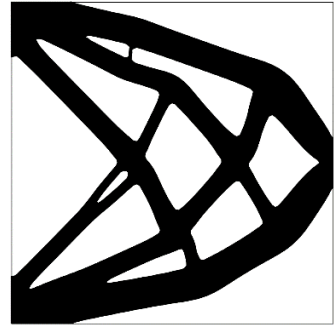 |

Figure 19. Diverse and competitive designs generated by DNN-based level set method

## 5. Conclusion

In this paper, a DNN level set method is proposed for topology optimization. Deep neural networks are popular for function approximation. The implicit function is represented by deep feedforward neural networks. The activation function is chosen as hyperbolic tangent function. Based on deep neural networks, a high level of smoothness in gradient and curvature of implicit function can be achieved. The Hamilton-Jacobi PDE is transformed into parametrized ODE and implicit function is updated through updating the weights and biases of network.

The main contribution of the proposed method is applying the DNN as a function approximator to describe implicit function of the level set method. Different DNN architectures are capable of generating diverse and competitive designs with excellent structural performance. DNN-based level set method can provide the designers with multiple conceptual alternatives instead of finding one optimum solution to maximize or minimize the objective. The limitation of the present work is that the mathematical connection between network architecture and structural complexity or performance cannot be quantified explicitly. At present, applying mathematical tool to quantify this relation is extremely difficult and further effort will be devoted to solving this issue in the future. In addition, using deep learning method to represent implicit function opens an opportunity toward a marriage of machine learning and topology optimization.

**Acknowledgement**

Financial support from the National Science Foundation (CMMI-1634261) for this work is gratefully acknowledged.

# Reference

[1] M. P. Bendsoe and O. Sigmund, *Topology optimization: theory, methods, and applications*. Springer Science & Business Media, 2013.

[2] O. Sigmund and K. Maute, "Topology optimization approaches," *Structural and Multidisciplinary Optimization,* vol. 48, no. 6, pp. 1031-1055, 2013.

[3] M. P. Bendsøe and O. Sigmund, "Material interpolation schemes in topology optimization," *Archive of applied mechanics,* vol. 69, no. 9-10, pp. 635-654, 1999.

[4] N. P. van Dijk, K. Maute, M. Langelaar, and F. Van Keulen, "Level-set methods for structural topology optimization: a review," *Structural and Multidisciplinary Optimization,* vol. 48, no. 3, pp. 437-472, 2013.

[5] F. Wang, B. S. Lazarov, and O. Sigmund, "On projection methods, convergence and robust formulations in topology optimization," *Structural and Multidisciplinary Optimization,* vol. 43, no. 6, pp. 767-784, 2011.

[6] J. Norato, B. Bell, and D. A. Tortorelli, "A geometry projection method for continuum-based topology optimization with discrete elements," *Computer Methods in Applied Mechanics and Engineering,* vol. 293, pp. 306-327, 2015.

[7] S. Watts and D. A. Tortorelli, "A geometric projection method for designing three‐dimensional open lattices with inverse homogenization," *International Journal for Numerical Methods in Engineering,* vol. 112, no. 11, pp. 1564-1588, 2017.

[8] B. S. Lazarov and F. Wang, "Maximum length scale in density based topology optimization," *Computer Methods in Applied Mechanics and Engineering,* vol. 318, pp. 826-844, 2017.

[9] M. Zhou, B. S. Lazarov, F. Wang, and O. Sigmund, "Minimum length scale in topology optimization by geometric constraints," *Computer Methods in Applied Mechanics and Engineering,* vol. 293, pp. 266-282, 2015.

[10] B. S. Lazarov, F. Wang, and O. Sigmund, "Length scale and manufacturability in density-based topology optimization," *Archive of Applied Mechanics,* vol. 86, no. 1-2, pp. 189-218, 2016.

[11] B. S. Lazarov, M. Schevenels, and O. Sigmund, "Robust design of large-displacement compliant mechanisms," *Mechanical sciences,* vol. 2, no. 2, pp. 175-182, 2011.

[12] J. K. Guest, "Topology optimization with multiple phase projection," *Computer Methods in Applied Mechanics and Engineering,* vol. 199, no. 1-4, pp. 123-135, 2009.

[13] J. K. Guest, "Imposing maximum length scale in topology optimization," *Structural and Multidisciplinary Optimization,* vol. 37, no. 5, pp. 463-473, 2009.

[14] J. K. Guest, J. H. Prévost, and T. Belytschko, "Achieving minimum length scale in topology optimization using nodal design variables and projection functions," *International journal for numerical methods in engineering,* vol. 61, no. 2, pp. 238-254, 2004.

[15] A. Asadpoure, M. Tootkaboni, and J. K. Guest, "Robust topology optimization of structures with uncertainties in stiffness–Application to truss structures," *Computers & Structures,* vol. 89, no. 11-12, pp. 1131-1141, 2011.

[16] M. Schevenels, B. S. Lazarov, and O. Sigmund, "Robust topology optimization accounting for spatially varying manufacturing errors," *Computer Methods in Applied Mechanics and Engineering,* vol. 200, no. 49-52, pp. 3613-3627, 2011.

[17] O. Sigmund, "Manufacturing tolerant topology optimization," *Acta Mechanica Sinica,* vol. 25, no. 2, pp. 227-239, 2009.

[18] B. S. Lazarov, M. Schevenels, and O. Sigmund, "Topology optimization with geometric uncertainties by perturbation techniques," *International Journal for Numerical Methods in Engineering,* vol. 90, no. 11, pp. 1321-1336, 2012.

[19] O. Sigmund, "Morphology-based black and white filters for topology optimization," *Structural and Multidisciplinary Optimization,* vol. 33, no. 4-5, pp. 401-424, 2007.

[20] J. A. Sethian, "Theory, algorithms, and applications of level set methods for propagating interfaces," *Acta numerica,* vol. 5, pp. 309-395, 1996.
[21] S. Osher and J. A. Sethian, "Fronts propagating with curvature-dependent speed: algorithms based on Hamilton-Jacobi formulations," *Journal of computational physics,* vol. 79, no. 1, pp. 12-49, 1988.
[22] S. J. Osher and F. Santosa, "Level set methods for optimization problems involving geometry and constraints: I. Frequencies of a two-density inhomogeneous drum," *Journal of Computational Physics,* vol. 171, no. 1, pp. 272-288, 2001.
[23] G. Allaire, F. Jouve, and A.-M. Toader, "Structural optimization using sensitivity analysis and a level-set method," *Journal of computational physics,* vol. 194, no. 1, pp. 363-393, 2004.
[24] M. Y. Wang, X. Wang, and D. Guo, "A level set method for structural topology optimization," *Computer methods in applied mechanics and engineering,* vol. 192, no. 1-2, pp. 227-246, 2003.
[25] S. Wang and M. Y. Wang, "Radial basis functions and level set method for structural topology optimization," *International journal for numerical methods in engineering,* vol. 65, no. 12, pp. 2060-2090, 2006.
[26] P. Wei and M. Y. Wang, "Piecewise constant level set method for structural topology optimization," *International Journal for Numerical Methods in Engineering,* vol. 78, no. 4, pp. 379-402, 2009.
[27] L. Jiang, S. Chen, and X. Jiao, "Parametric shape and topology optimization: A new level set approach based on cardinal basis functions," *International Journal for Numerical Methods in Engineering,* vol. 114, no. 1, pp. 66-87, 2018.
[28] X. Guo, W. Zhang, and W. Zhong, "Doing topology optimization explicitly and geometrically—a new moving morphable components based framework," *Journal of Applied Mechanics,* vol. 81, no. 8, 2014.
[29] W. Zhang, J. Zhou, Y. Zhu, and X. Guo, "Structural complexity control in topology optimization via moving morphable component (MMC) approach," *Structural and Multidisciplinary Optimization,* vol. 56, no. 3, pp. 535-552, 2017.
[30] W. Zhang *et al.*, "Explicit three dimensional topology optimization via Moving Morphable Void (MMV) approach," *Computer Methods in Applied Mechanics and Engineering,* vol. 322, pp. 590-614, 2017.
[31] W. Zhang, D. Li, J. Zhou, Z. Du, B. Li, and X. Guo, "A moving morphable void (MMV)-based explicit approach for topology optimization considering stress constraints," *Computer Methods in Applied Mechanics and Engineering,* vol. 334, pp. 381-413, 2018.
[32] L. Jiang and S. Chen, "Parametric structural shape & topology optimization with a variational distance-regularized level set method," *Computer Methods in Applied Mechanics and Engineering,* vol. 321, pp. 316-336, 2017.
[33] Y. Luo, J. Xing, and Z. Kang, "Topology optimization using material-field series expansion and Kriging-based algorithm: An effective non-gradient method," *Computer Methods in Applied Mechanics and Engineering,* vol. 364, p. 112966, 2020.
[34] P. Lison, "An introduction to machine learning," *Language Technology Group: Edinburgh, UK,* 2015.
[35] M. Rastegari, V. Ordonez, J. Redmon, and A. Farhadi, "Xnor-net: Imagenet classification using binary convolutional neural networks," in *European conference on computer vision*, 2016: Springer, pp. 525-542.
[36] E. Gawehn, J. A. Hiss, and G. Schneider, "Deep learning in drug discovery," *Molecular informatics,* vol. 35, no. 1, pp. 3-14, 2016.
[37] M. Raissi, P. Perdikaris, and G. E. Karniadakis, "Physics-informed neural networks: A deep learning framework for solving forward and inverse problems involving nonlinear partial differential equations," *Journal of Computational Physics,* vol. 378, pp. 686-707, 2019.
[38] M. Raissi, A. Yazdani, and G. E. Karniadakis, "Hidden fluid mechanics: Learning velocity and pressure fields from flow visualizations," *Science,* vol. 367, no. 6481, pp. 1026-1030, 2020.

[39] R. Iten, T. Metger, H. Wilming, L. Del Rio, and R. Renner, "Discovering physical concepts with neural networks," *Physical Review Letters,* vol. 124, no. 1, p. 010508, 2020.
[40] S. L. Brunton, B. R. Noack, and P. Koumoutsakos, "Machine learning for fluid mechanics," *Annual Review of Fluid Mechanics,* vol. 52, pp. 477-508, 2020.
[41] M. Raissi, Z. Wang, M. S. Triantafyllou, and G. E. Karniadakis, "Deep learning of vortex-induced vibrations," *Journal of Fluid Mechanics,* vol. 861, pp. 119-137, 2019.
[42] Y. Yu, T. Hur, J. Jung, and I. G. Jang, "Deep learning for determining a near-optimal topological design without any iteration," *Structural and Multidisciplinary Optimization,* vol. 59, no. 3, pp. 787-799, 2019.
[43] X. Lei, C. Liu, Z. Du, W. Zhang, and X. Guo, "Machine learning-driven real-time topology optimization under moving morphable component-based framework," *Journal of Applied Mechanics,* vol. 86, no. 1, 2019.
[44] S. Oh, Y. Jung, S. Kim, I. Lee, and N. Kang, "Deep generative design: Integration of topology optimization and generative models," *Journal of Mechanical Design,* vol. 141, no. 11, 2019.
[45] J. J. Park, P. Florence, J. Straub, R. Newcombe, and S. Lovegrove, "Deepsdf: Learning continuous signed distance functions for shape representation," in *Proceedings of the IEEE Conference on Computer Vision and Pattern Recognition*, 2019, pp. 165-174.
[46] B. Wang, Y. Zhou, Y. Zhou, S. Xu, and B. Niu, "Diverse competitive design for topology optimization," *Structural and Multidisciplinary Optimization,* vol. 57, no. 2, pp. 891-902, 2018.
[47] K. Yang *et al.*, "Simple and effective strategies for achieving diverse and competitive structural designs," *Extreme Mechanics Letters,* vol. 30, p. 100481, 2019.
[48] Y. He, K. Cai, Z.-L. Zhao, and Y. M. Xie, "Stochastic approaches to generating diverse and competitive structural designs in topology optimization," *Finite Elements in Analysis and Design,* vol. 173, p. 103399, 2020.
[49] I. Goodfellow, Y. Bengio, and A. Courville, *Deep learning*. MIT press, 2016.
[50] G. Cybenko, "Approximation by superpositions of a sigmoidal function," *Mathematics of control, signals and systems,* vol. 2, no. 4, pp. 303-314, 1989.
[51] A. Paszke *et al.*, "Automatic differentiation in pytorch," 2017.
[52] D. E. Rumelhart, G. E. Hinton, and R. J. Williams, "Learning internal representations by error propagation," California Univ San Diego La Jolla Inst for Cognitive Science, 1985.
[53] J. C. Butcher, *The numerical analysis of ordinary differential equations: Runge-Kutta and general linear methods*. Wiley-Interscience, 1987.
[54] S. Wang, K. M. Lim, B. C. Khoo, and M. Y. Wang, "An extended level set method for shape and topology optimization," *Journal of Computational Physics,* vol. 221, no. 1, pp. 395-421, 2007.
[55] S. Osher, R. Fedkiw, and K. Piechor, "Level set methods and dynamic implicit surfaces," *Appl. Mech. Rev.,* vol. 57, no. 3, pp. B15-B15, 2004.
[56] M. Y. Wang and P. Wang, "The augmented Lagrangian method in structural shape and topology optimization with RBF based level set method," in *CJK-OSM 4: The Fourth China-Japan-Korea Joint Symposium on Optimization of Structural and Mechanical Systems*, 2006, p. 191.
[57] C. Li, C. Xu, C. Gui, and M. D. Fox, "Distance regularized level set evolution and its application to image segmentation," *IEEE transactions on image processing,* vol. 19, no. 12, pp. 3243-3254, 2010.
[58] V. J. Challis, "A discrete level-set topology optimization code written in Matlab," *Structural and multidisciplinary optimization,* vol. 41, no. 3, pp. 453-464, 2010.
[59] D. Hartmann, M. Meinke, and W. Schröder, "The constrained reinitialization equation for level set methods," *Journal of computational physics,* vol. 229, no. 5, pp. 1514-1535, 2010.
[60] G. A. Anastassiou, "Multivariate hyperbolic tangent neural network approximation," *Computers & Mathematics with Applications,* vol. 61, no. 4, pp. 809-821, 2011.
[61] E. Andreassen, A. Clausen, M. Schevenels, B. S. Lazarov, and O. Sigmund, "Efficient topology optimization in MATLAB using 88 lines of code," *Structural and Multidisciplinary Optimization,* vol. 43, no. 1, pp. 1-16, 2011.